\documentclass[a4paper,12pt]{article}
\usepackage{amsmath,amssymb,graphicx}
\usepackage{amsthm}
\usepackage{algorithm}
\usepackage{algorithmic}
\usepackage{caption,color}
\usepackage{hyperref}
\usepackage{tikz}
\usepackage{subcaption}

\allowdisplaybreaks[3]

\title{Spectral extremal graphs for $F_6$-free graphs with even size\thanks{Email:
\url{Ylj99@hnu.edu.cn} (L. Yu),
,\url{ypeng1@hnu.edu.cn} (Y. Peng, corresponding author)}}
\author{\small{Loujun Yu$^1$,  Yuejian Peng$^1$}\\
{\small $^1$School of Mathematics, Hunan University}\\
{\small Changsha, Hunan, 410082, P.R. China} \\
\makeatletter
}
\newcommand\tabcaption{\def\@captype{table}\caption}
\date{\today}
\makeatother

\newtheorem{theorem}{Theorem}[section]

\newtheorem*{rmk}{Remark}
\newtheorem{lemma}[theorem]{Lemma}

\newtheorem{prob}[theorem]{Problem}

\newtheorem{conjecture}[theorem]{Conjecture}

\newtheorem{claim}{Claim}

\newcommand*{\QED}{\hfill\ensuremath{\square}}

\begin{document}
\maketitle
\begin{abstract}
Let $F_l$ be the fan graph obtained by joining a vertex with a path on $l-1$ vertices. Yu, Li and Peng [Discrete Math. 346 (2023)]  conjectured that if the number of edges  of $G$ is $m$ and the spectral radius $\lambda(G)>\frac{k-1+\sqrt{4m-k^2+1}}{2}$, then $G$ contains a $F_{2k+1}$ and $F_{2k+2}$, unless $G=K_{k}\vee (\frac{m}{k}-\frac{k-1}{2})K_1$.
The case $k\geq 3$ of the above conjecture has been confirmed by Li, Zhao and Zou [J. Graph theory 110 (2025)]. Zhang and Wang [Discrete Math. 347 (2024)], Yu, Li and Peng [Discrete Math. 348 (2025)], Gao and Li [Discrete Math. 349 (2026)] confirmed the  case $k=2$. However, the extremal graphs for the case $k=2$ only exist when $m$ is odd. The case with $m$ even  has not been determined.
In this paper, we characterize the extremal graph for $F_6$ and even  $m\ge 3000$.
\end{abstract}



\section{Introduction}

As usual, let $K_n$, $C_n$ and $P_n$ be the complete graph, the cycle and the path on $n$ vertices, respectively.
Let $K_{s,t}$ be the complete bipartite graph with parts of sizes
$s$ and $t$. In particular, $K_{1,r}$ is called a star and the vertex  with degree $r$ is called the center vertex.
A double star is a graph obtained by taking  an edge and joining one of its end vertices with $a$ new vertices and the other end vertex with $b$ new vertices, where the vertices with degrees $a+1$ and $b+1$ are called the center vertices of the double star.
Let $K_n-e$  and $K_{1,r}+e$ denote the graphs obtained by deleting an edge from $K_n$ and by adding an edge within two independent vertices of $K_{1,r}$, respectively.
Let $kG$ denote the union
of $k$ vertex-disjoint copies of $G$.
Let $G \vee H$ be the join graph obtained by taking the vertex disjoint union of $G$ and $H$ and joining  each vertex of $G$
 to each vertex of $H$.

For an $n$-vertex graph $G$, the adjacency matrix  is defined as $A(G)=[a_{i,j}]_{i,j=1}^{n}$, where $a_{i,j}=1$ if $\{i,j\}\in E(G)$, and $a_{i,j}=0$ otherwise. The spectral radius, denoted by $\lambda(G)$, is the maximum modulus of eigenvalues of $A(G)$. From the definition of adjacency matrix, it is obvious that $A(G)$ is non-negative. By the Perron-Frobenius theorem, $\lambda(G)$ is actually the largest eigenvalue of $A(G)$ and has a corresponding non-negative eigenvector. Moreover, if $G$ is connected, then $\lambda(G)$ has a unique positive eigenvector up to a positive scaling coefficient, which is called the Perron vector.

Dating back to 1970,  Nosal \cite{Nosal1970} (see, e.g., \cite{Ning2017-ars} for alternative proofs) showed that for
every triangle-free graph $G$ with $m$ edges, we have
$\lambda (G)\le \sqrt{m} $.
Later, Nikiforov \cite{Niki2002cpc,Niki2006-walks,Niki2009jctb} extended Nosal's result by
proving that
if $G$ is $K_{r+1}$-free, then
\begin{equation}       \label{eq1}
\lambda (G) \le \sqrt{2m\Bigl( 1-\frac{1}{r}\Bigr)}.
\end{equation}
Moreover, the equality holds
if and only if
$G$ is a complete bipartite graph for $r=2$,
or a regular complete $r$-partite graph for each $r\ge 3$.
Tur\'{a}n's theorem
can be derived from  (\ref{eq1}).
Indeed, using Rayleigh's inequality, we have
$\frac{2m}{n}\le \lambda (G)\le  (1-\frac{1}{r})^{1/2} \sqrt{2m}$,
which yields $ m \le (1- \frac{1}{r}) \frac{n^2}{2}$.
Thus,  (\ref{eq1}) could be viewed as
a spectral extension of Tur\'{a}n's theorem.
Furthermore, Bollob\'{a}s and Nikiforov \cite{BN2007jctb}
 conjectured that
if  $G$ is a $K_{r+1}$-free graph
of order at least $r+1$
with $m$ edges, then
\begin{equation} \label{eq-BN}
{  \lambda_1^2(G)+ \lambda_2^2(G)
\le 2m\Bigl( 1-\frac{1}{r}\Bigr) }.
\end{equation}
This conjecture has attracted wide attention and has been confirmed for some special cases, see
 \cite{AL2015,LSY2022,LNW2021,Niki2021,Zhang2024}.
Nevertheless this intriguing
 problem remains open,
we refer the readers to
\cite{ELW2024}
for two related conjectures.
Both (\ref{eq1}) and (\ref{eq-BN})
boosted the great interests of studying the Brualdi-Hoffman-Tur\'{a}n problem, which is stated as follows.
\begin{prob}[Brualdi-Hoffman-Tur\'{a}n problem] Given a graph $F$ and a positive integer $m$,
what is the maximum spectral radius of an $F$-free graph with $m$ edges?
\end{prob}
The graph attaining maximum spectral radius among all $F$-free graphs with $m$ edges is called an extremal graph for $F$.

This problem has been studied for various families of graphs. For example, see \cite{LFP2023-solution,LP2022second,LLL2022,MLH2022,Niki2009laa,Wang2022DM,WG2024,ZLS2021,ZS2022dm} for  $C_k$-free graphs,
see \cite{ZLS2021} for  $K_{2,r+1}$-free graphs, see \cite{FYH2022,LZS2024,LZZ2025JGT,LW2022,LSW2022} for $C_{k}^+$-free graphs,
where $C_k^{+}$
is a graph on $k$ vertices obtained from $C_k$
by adding a chord between two vertices with distance two;
see \cite{Niki2021,NZ2021} for  $B_k$-free graphs,
where $B_k$ denotes the book graph which consists  of $k$ triangles
sharing a common edge,
see \cite{CY2025DM,GL2026,LZZ2025JGT,YLP2025DM,ZW2024DM,ZP2022} for fan graphs,
see \cite{LLP2022,LZZ2025JGT,YLP2025DM} for $F_{k,3}$-free graphs,
where $F_{k,3}$
is the friendship graph  consisting of $k$ triangles  intersecting in
a common vertex, see \cite{GL2025,GL2025-2,LZZ2025JGT} for theta graphs, where a theta graph, say $\theta_{r,p,q}$, is the graph obtained by connecting two distinct vertices  with three internally disjoint paths of length $r$, $p$, $q$, where $q\geq p\geq r\geq 1$ and $p\geq 2$.

In particular, we aim to consider the Brualdi-Hoffman-Tur\'{a}n problem of a fan graph, denoted by $F_l$, which is obtained by joining an isolated vertex to $P_{l-1}$.
 In 2025, Yu, Li and Peng \cite{YLP2025DM} proposed the following  conjecture.

\begin{conjecture}[Yu--Li--Peng \cite{YLP2025DM}, 2025]\label{con-YLP}
Let $k\geq 2$ be fixed and $m$ be sufficiently large.
If $G$ is an $F_{2k+1}$-free or $F_{2k+2}$-free graph with $m$ edges without isolated vertices, then
\[ \lambda (G)\le \frac{k-1+\sqrt{4m-k^2+1}}{2}, \]
 and the equality holds if and only if $G = K_k\vee (\frac{m}{k}-\frac{k-1}{2})K_1$.
\end{conjecture}

Subsequently, Zhang and Wang \cite{ZW2024DM}, Yu, Li and Peng \cite{YLP2025DM} confirmed the conjecture for $F_5$. And Li, Zhao and Zou \cite{LZZ2025JGT} solved the conjecture for $k\geq 3$.

\begin{theorem}[Zhang--Wang, Yu--Li--Peng \cite{YLP2025DM,ZW2024DM}, 2024, 2025]\label{Thm-ZWYLP}
Let $G$ be an $F_5$-free graph with $m\geq 11$ edges without isolated vertices. Then
\[ \lambda (G)\le \frac{1+\sqrt{4m-3}}{2}, \]
 and the equality holds if and only if $G = K_2\vee \frac{m-1}{2}K_1$.
\end{theorem}

\begin{theorem}[Li--Zhao--Zou \cite{LZZ2025JGT}, 2025]\label{LZZ25}
Let $k\geq 3$ and $m\geq \frac{9}{4}k^6+6k^5+46k^4+56k^3+196k^2$.
If $G$ is an $F_{2k+1}$-free or $F_{2k+2}$-free graph with $m$ edges without isolated vertices, then
\[ \lambda (G)\le \frac{k-1+\sqrt{4m-k^2+1}}{2}, \]
 and the equality holds if and only if $G = K_k\vee (\frac{m}{k}-\frac{k-1}{2})K_1$.
\end{theorem}

\begin{rmk}
Obviously, for $k\geq 3$, the extremal graph $K_k\vee (\frac{m}{k}-\frac{k-1}{2})K_1$ only exists when $\frac{m}{k}-\frac{k-1}{2}$ is an integer. If  $\frac{m}{k}-\frac{k-1}{2}$ is not an integer,  extremal graphs for $F_{2k+1}$ and $F_{2k+2}$ are still not known.
\end{rmk}

Recently, Gao and Li \cite{GL2026} proved the remaining case $F_6$ of Conjecture \ref{con-YLP} .

\begin{theorem}[Gao--Li \cite{GL2026}, 2026]\label{Thm-GL}
Let $G$ be an $F_6$-free graph with $m\geq 88$ edges without isolated vertices. Then
\[ \lambda (G)\le \frac{1+\sqrt{4m-3}}{2}, \]
 and the equality holds if and only if $G = K_2\vee \frac{m-1}{2}K_1$.
\end{theorem}

However, it is easy to see that the extremal graphs in Theorems \ref{Thm-ZWYLP} and \ref{Thm-GL} only exist when $m$ is odd. The case of even size has not been determined. Let $S_{n,k}$  be the graph obtained by joining $K_k$ to $n-k$ isolated vertices. Obviously, $K_2\vee \frac{m-1}{2}K_1=S_{\frac{m+3}{2},2}$. Let $S_{n,2}^{-t}$ be the graph obtained from $S_{n,2}$ by deleting $t$ edges between a vertex of degree $n-1$ and $t$ vertices of degree 2. Recently, Chen and Yuan \cite{CY2025DM} characterized the extremal graph for $F_5$ with given even size.

\begin{theorem}[Chen--Yuan \cite{CY2025DM}, 2025]\label{le-Chen}
Let $G$ be an $F_5$-free graph with $m\geq 92$ edges without isolated vertices. If $m$ is even, then
\[ \lambda (G)\le \lambda(S_{\frac{m+4}{2},2}^{-1}), \]
 and the equality holds if and only if $G = S_{\frac{m+4}{2},2}^{-1}$.
\end{theorem}

In this note, we characterize the extremal graph for $F_6$ with given even size.

\begin{theorem}\label{main-Thm}
If $G$ is an $F_6$-free graph with $m$ edges without isolated vertices, then for even $m\geq 3000$,
\[ \lambda (G)\le \lambda(S_{\frac{m+4}{2},2}^{-1}), \]
 and the equality holds if and only if $G = S_{\frac{m+4}{2},2}^{-1}$.
\end{theorem}

\section{Preliminaries}

 Given a vertex $u\in V(G)$,
 {let $N_G(u)$ be the set of neighbors of $u$,
and $d_G(u)$ be the degree of $u$ in $G$.
Moreover, we denote  $N_G[u]=N_G(u)\cup\{u\}$.
Sometimes, we will eliminate the subscript and write $N(u)$ and $d(u)$ if there is no confusion.
For a subset $U\subseteq V(G)$, we write $e(U)$ for the
number of edges with two endpoints in $U$.
For two disjoint sets $U,W$, we write $e(U,W)$ for the number of edges between $U$ and $W$.
For simplicity,  we denote by
$N_{U}(u)$ the set of vertices in $U$ that are  adjacent to $u$, i.e, $N_U(u)=N_G(u)\cap U$, and let $d_U(u)$ be the number of vertices in $N_U(u)$.

 \begin{lemma}[Wu--Xiao--Hong \cite{WXH2005}, 2005] \label{le-niki}
Let $G$ be a connected graph
and $(x_1,\ldots ,x_n)^\mathrm{T}$ be a Perron vector of $G$,
where the coordinate $x_i$ corresponds to the vertex $v_i$.
Assume that
 $v_i,v_j \in V(G)$ are vertices such that $x_i \ge x_j$, and $S\subseteq N_G[v_j] \setminus N_G[v_i]$ is non-empty.
 Denote $G'=G- \{v_jv : v\in S\} +
\{v_iv : v\in S\}$. Then $\lambda (G) < \lambda (G')$.
\end{lemma}

 An extremal vertex is a vertex corresponding {to} the maximum coordinate of the Perron vector of the extremal graph.
\begin{lemma}[Zhai--Lin--Shu \cite{ZLS2021}, 2021]\label{le-zhai}
Let $F$ be a $2$-connected graph and
$G^*$ be an extremal graph for $F$ with an extremal vertex $u^*$. Then $G^*$ is connected and there is no cut vertex in $V(G^*)\setminus\{u^*\}$. Furthermore,
we have $d(u)\ge  2$ for any vertex $u\in V(G^*)\setminus N[u^*]$.
\end{lemma}

\begin{lemma}[Fang--You \cite{FYH2022}, 2021]\label{le-Fang}
For $m>t+1$,
\begin{align*}
\lambda (S_{\frac{m+t+2}{2},2}^{-t})<
\left\{
\begin{array}{ll}
\lambda (S_{\frac{m+3}{2},2})& \text{if $t (\geq 2) $ is even,}\\
\lambda (S_{\frac{m+4}{2},2}^{-1}) & \text{if $t (\geq 3) $ is odd.}
\end{array}
\right.
\end{align*}
\end{lemma}

\section{Proof of Theorem \ref{main-Thm}}
  Let $m$ be even and let $G^*$ denote an extremal graph for $F_6$ with $m$ edges. For convenience, we will use {$\lambda$} and $\mathbf{x}$ to denote the spectral radius and the Perron vector of $G^*$, where the coordinate $x_u$ corresponds to the vertex $u$.
Let $u^*$ be an extremal vertex of $G^*$,
 $U=N(u^*)$ and $W=V(G^*)\setminus N[u^*]$.  Since $F_6$ is 2-connected, by Lemma \ref{le-zhai}, we know that $G^*$ is connected. In addition,
we get $\lambda \ge  \lambda (S_{\frac{m+4}{2},2}^{-1})>\frac{1+\sqrt{4m-5}}{2}>55$ since $S_{\frac{m+4}{2},2}^{-1}$ is $F_6$-free and $m\geq 3000$.
 Hence, $\lambda^2-\lambda  >  m-\frac{3}{2}$.
 Since
 \[ \lambda^2 x_{u^*} = |U| x_{u^*} +
\sum_{u\in U} d_U(u)x_u + \sum_{w\in W} d_U(w) x_w, \]
which together with $\lambda x_{u^*}=\sum_{u\in U}x_u$  yields
\[  (\lambda^2 - \lambda ) x_{u^*}
= |U| x_{u^*} + \sum_{u\in U} (d_U(u) -1)x_u
 + \sum_{w\in W} d_U(w) x_w.  \]
Note that $\lambda^2 - \lambda > m-\frac{3}{2}
= |U| + e(U) + e(U,W) +e(W) -\frac{3}{2}$. Then
  \begin{align}\label{eqc}
   \left(e(U)+e(W)+e(U,W)-\frac{3}{2}\right)x_{u^*} <  \sum\limits_{u\in U}(d_U(u)-1)x_u+\sum\limits_{w\in W}d_U(w)x_w.
  \end{align}
By simplifying, we can get
\begin{equation} \label{eqa}
 e(W)< \sum\limits_{u\in U}(d_U(u)-1)\frac{x_u}{x_{u^*}}
 -  e(U) + \sum_{w\in W} d_U(w) \frac{x_w}{x_{u^*}} - e(U,W) +\frac{3}{2}.
\end{equation}

\begin{lemma} \label{le-eU}
$e(U)\ge 1$.
\end{lemma}

\noindent \textbf{Proof}.
As $\lambda \ge \frac{1+\sqrt{4m-5}}{2}> \sqrt{m}$ when $m\geq 2$, we have
\begin{align*}
mx_{u^*}& < \lambda^{2}x_{u^*}=|U|x_{u^*}+\sum\limits_{u\in U}d_U(u)x_u+\sum\limits_{w\in W}d_U(w)x_w \\
&\le |U| x_{u^*} + 2e(U) x_{u^*} + e(U,W)x_{u^*}.
\end{align*}
Observe that $m=|U| + e(U) + e(U,W) + e(W)$.
Therefore, the above inequality  implies that $e(W)<e(U)$. If $e(U)=0$, then $e(W)<0$, a contradiction. Thus, $e(U)\ge  1$. \QED

\medskip

An isolated vertex in $G^*[U]$ is called a trivial component.
Let $U_0$ and $W_0$ denote the set of isolated vertices in $G^*[U]$ and $G^*[W]$, respectively.
By Lemma \ref{le-eU},
there {exists} at least one non-trivial  {component} in $G^*[U]$.
Let $U_+=\{u\in U|~d_U(u)>0\}$.
Then (\ref{eqa}) can be written as
\begin{equation} \label{eq-etb}
 e(W)< \sum\limits_{u\in U_+}(d_{U_+}(u)-1)\frac{x_u}{x_{u^*}}
 - e(U_+)-\sum_{u\in U_0}\frac{x_u}{x_{u^*}} + \sum_{w\in W} d_U(w) \frac{x_w}{x_{u^*}} - e(U,W) +\frac{3}{2}.
\end{equation}
Let $\mathcal{H}$ be the set of all components in $G^*[U_+]$.  For each non-trivial component $H\in \mathcal{H}$, we define $W_H=\cup_{u\in V(H)}N_W(u)$ and
\begin{equation} \label{eq-eta}
 \eta(H) :=\sum\limits_{u\in V(H)}(d_H(u)-1)\frac{x_u}{x_{u^*}}-e(H).
 \end{equation}
 Clearly,  we have
$\eta (H) \le e(H) - |V(H)| $.
We define $\gamma (w)=d_U(w)(1-\frac{x_{w}}{x_{u^*}})$ for any $w\in W$ and $\gamma(W')=\sum_{w\in W'}\gamma(w)$ for any non-empty subset $W'\subseteq W$. Clearly, $\gamma(w)$ and $\gamma(W)$ are nonnegative.
Then (\ref{eq-etb}) can be rewritten as
\begin{align}\label{eqb}
 e(W) <  \sum\limits_{H\in \mathcal{H}}\eta(H)-\sum\limits_{u\in U_0}\frac{x_u}{x_{u^*}}-\gamma(W)+\frac{3}{2}.
\end{align}
Recall that $G^*$ is $F_6$-free, we  know that $G^*[U]$ is $P_5$-free.
 It implies that any component of $\mathcal{H}$ is isomorphic to  a star, or a double star, or $K_{1,r}+e$ for some integer $r$, or $C_4$, or $K_4-e$, or $K_4$.  Since  $\eta (H) \le e(H) - |V(H)| $ holds for any $H\in \mathcal{H}$,  we immediately get the following lemma by simple calculations.
 \begin{lemma}\label{le-etah}
 For any $H\in \mathcal{H}$,
 \begin{align*}
 \eta(H)\leq \left\{
 \begin{array}{ll}
 -1, & \text{if $H$ is a star or a double star};\\
 0, & \text{if $H=K_{1,r}+e$ for some integer $r$ or a $C_4$};\\
 1, & \text{if $H=K_{4}-e$};\\
 2, & \text{if $H=K_4$}.\\
 \end{array}
 \right.
 \end{align*}
 \end{lemma}

 \begin{lemma}\label{le-etah0}
 $\eta(H)\leq 0$ holds for any $H\in \mathcal{H}$.
 \end{lemma}
 \noindent \textbf{Proof}.  Let $\mathcal{H}_+=\{H\in \mathcal{H}| ~\eta(H)>0\}$. By Lemma \ref{le-etah}, we know that $H= K_4-e$ or $K_4$ for any $H\in \mathcal{H}_+$.

\begin{claim}\label{cla-etah0}
$\sum_{u\in V(H)}x_u\geq\frac{5}{2}x_{u^*}$ for any $H\in \mathcal{H}_+$.
\end{claim}

\noindent{Proof of Claim \ref{cla-etah0}}. Suppose to the contrary that there is a component $H\in \mathcal{H}_+$ such that $\sum_{u\in V(H)}x_u<\frac{5}{2}x_{u^*}$.  Since $H= K_4-e$ or $K_4$,  we have $\Delta(H)\leq 3$ and $e(H)\geq 5$. Furtheremore,
 \begin{align*}
 \eta(H)&=\sum\limits_{u\in V(H)}(d_H(u)-1)\frac{x_u}{x_{u^*}}-e(H)\leq 2 \sum\limits_{u\in V(H)}\frac{x_u}{x_{u^*}}-5<0,
 \end{align*}
 a contradiction to the definition of $\mathcal{H}_+$.\QED
\begin{claim}\label{cla-h1}
$|\mathcal{H}_{+}|< \frac{2\lambda^2-2\lambda +3}{5\lambda-11}$.
\end{claim}
\noindent {Proof of Claim \ref{cla-h1}}. For any $H\in \mathcal{H}_+$,  we will show that $e(H,W)\geq \frac{5\lambda-29}{2}$. Indeed, we could see that
\begin{align*}
\lambda\sum\limits_{u\in V(H)}x_u &=\sum\limits_{u\in V(H)}\left(x_{u^*}+\sum\limits_{v\in N_{H}(u)}x_v+\sum\limits_{w\in N_{W}(u)}x_w\right)\\
&=4x_{u^*}+\sum\limits_{u\in V(H)}\sum\limits_{v\in N_{H}(u)}x_v+\sum\limits_{u\in V(H)}\sum\limits_{w\in N_{W}(u)}x_w\\
&\leq 4x_{u^*}+\sum\limits_{u\in V(H)}d_H(u)x_u+e(H,W)x_{u^*}\\
&=4x_{u^*}+\eta(H) x_{u^*}+e(H) x_{u^*}+\sum\limits_{u\in V(H)}x_u+e(H,W)x_{u^*}.
\end{align*}
Hence, by Lemma \ref{le-etah} and $e(H)\leq 6$, we have
\begin{align}\label{eq-hwe}
(\lambda-1) \sum\limits_{u\in V(H)}x_u&\leq 4x_{u^*}+\eta(H) x_{u^*}+e(H) x_{u^*}+e(H,W)x_{u^*}\nonumber\\
&\leq 12 x_{u^*}+e(H,W) x_{u^*}.
\end{align}
Combining (\ref{eq-hwe}) with Claim \ref{cla-etah0}, we have $e(H,W)\geq \frac{5\lambda-29}{2}$. Thus,
\begin{align*}
\lambda^2-\lambda+\frac{3}{2}>m&\geq \sum\limits_{H\in \mathcal{H}_+} \left(|H|+e(H)+e(H,W)\right)\\
&\geq \left(4+5+\frac{5\lambda-29}{2}\right)|\mathcal{H}_+|\\
&=\frac{5\lambda-11}{2}|\mathcal{H}_+|,
\end{align*}
which implies that $|\mathcal{H}_+|<\frac{2\lambda^2-2\lambda +3}{5\lambda-11}$.\QED

If this lemma does not hold, then $\mathcal{H}_+\neq \varnothing$. Suppose that there exists a component $H'\in \mathcal{H}_+$ with $V(H')=\{u_1,u_2,u_3,u_4\}$. If $H'=K_4-e$, we will assume that $d_{H'}(u_2)=d_{H'}(u_3)=3$. Recall that $W_0=\{w\in W|~d_W(w)=0\}$. We will show that there is an upper bound for the entries of $\mathbf{x}$ corresponding to the vertices in $W_0$.

\begin{claim}\label{cla-w0}
If $W_0\neq \varnothing$, then $x_w\leq \frac{2\lambda-3}{2\lambda}x_{u^*}$ for any $w\in W_0$.
\end{claim}

\noindent {Proof of Claim \ref{cla-w0}}. For any $w\in W_0$, we have $N(w)\subseteq U=N(u^*)$. Suppose there exists a vertex $w_0\in W_0$ such that $x_{w_0}>\frac{2\lambda-3}{2\lambda}x_{u^*}$.

If $H'=K_4$, then $|N_{H'}(w_0)|\leq 1$ since $G^*$ is  $F_6$-free.  Combining it with Claim \ref{cla-etah0}, we know that $\lambda x_{w_0}\leq \lambda x_{u^*}-\frac{3}{2} x_{u^*}$, a contradiction to the choice of $w_0$.

If $H'=K_4-e$, then we claim that $w_0$ is adjacent to $u_2$ or $u_3$. Otherwise,  we have $\lambda x_{w_0}\leq \lambda x_{u^*}-x_{u_2}-x_{u_3}$. Since $x_{w_0}>\frac{2\lambda-3}{2\lambda}x_{u^*}$,   we will get $x_{u_2}+x_{u_3}<\frac{3}{2} x_{u^*}$. Moreover, $\eta(H')=\frac{2x_{u_1}+2x_{u_2}}{x_{u^*}}+\frac{x_{u_1}+x_{u_4}}{x_{u^*}}-5< 0$, which is a contradiction. Hence, $w_0$ is adjacent to $u_2$ or $u_3$. Without loss of generality, assume $w_0$ is adjacent to $u_2$, then  $N_{H'}(w_0)=\{u_2\}$ since $G^*$ is $F_6$-free. From Claim \ref{cla-etah0}, it is obvious that $\lambda x_{w_0}\leq \lambda x_{u^*}-x_{u_1}-x_{u_3}-x_{u_4}\leq \lambda x_{u^*}-\frac{3}{2} x_{u^*}$, a contradiction to the choice of $w_0$.\QED

Recall that  $W_{H'}=\cup_{i=1}^{4}N_{W}(u_i)$.
If $W_{H'}=\varnothing$, then $x_{u_i}\leq \frac{4}{\lambda}x_{u^*}$ for any $1\leq i\leq 4$. Then $\sum_{i=1}^{4}x_{u_i}\leq \frac{16}{\lambda}x_{u^*}$, which contradicts to Claim \ref{cla-etah0} as $\lambda>55$. Thus we have $W_{H'}\neq \varnothing$.

\begin{claim} \label{cla-enW4}
$e(W)\geq 4$.
\end{claim}

\noindent {Proof of Claim \ref{cla-enW4}}. Suppose to the contrary, $e(W)\leq 3$. We will prove the claim by distinguishing the following two cases.

 \textbf{Case 1}. $H'=K_4$.

Without loss of generality, suppose that $x_{u_1}\geq x_{u_2}\geq x_{u_3}\geq x_{u_4}$. Since $G^*$ is $F_6$-free, we have $N_W(u_i)\cap N_W(u_j)=\varnothing$ for any $1\leq i\neq j\leq 4$. Then we can get $N_W(u_i)=\varnothing$ for $2\leq i\leq 4$. Indeed, suppose there is  a vertex $w\in N_W(u_2)$. We can construct a graph $G'=G^*-u_2w+u_1w$. Clearly, $G'$ is an $F_6$-free graph with $m$ edges as $e(W)\leq 3$. And by Lemma \ref{le-niki}, $\lambda(G')>\lambda$, a contradiction. Furthermore, $x_{u_i}\leq \frac{4}{\lambda}x_{u^*}$ for $2\leq i\leq 4$. Then $\eta(H')=\sum_{i=1}^{4}\frac{2x_{u_i}}{x_{u^*}}-6\leq \frac{24}{\lambda}-4\leq 0$, which is contradicted to the choice of $H'$ as $\lambda>55$.

\textbf{Case 2}. $H'=K_4-e$.

Without loss of generality, suppose that $ x_{u_2}\geq x_{u_3}$. If $N_{W}(u_2)=\varnothing$ and  $N_{W}(u_3)=\varnothing$, then $x_{u_2}, x_{u_3}\leq \frac{4}{\lambda}x_{u^*}$. Thus, we obtain $\eta(H')=\frac{2x_{u_2}+2x_{u_3}}{x_{u^*}}+\frac{x_{u_1}+x_{u_4}}{x_{u^*}}-5\leq \frac{16}{\lambda}-3\leq 0$ as $\lambda>55$, a contradiction. Therefore, $N_{W}(u_2)\neq \varnothing $ or $N_{W}(u_3)\neq \varnothing$. In addition, we also have $N_{W}(u_2)\cap N_{W}(u_3)=\varnothing $ since $G^*$ is $F_6$-free. Similarly, we can get that $N_W(u_3)=\varnothing$. Otherwise, we can construct a new $m$-edge $F_6$-free graph $G'$ with  $\lambda(G')>\lambda$ as $x_{u_2}\geq x_{u_3}$ and $e(W)\leq 3$. Consequently, $x_{u_3}\leq \frac{4}{\lambda}x_{u^*}$ and $\eta(H')=\frac{2x_{u_2}+2x_{u_3}}{x_{u^*}}+\frac{x_{u_1}+x_{u_4}}{x_{u^*}}-5\leq \frac{8}{\lambda}-1\leq 0$ as $\lambda>55$, a contradiction to the choice of $H'$.\QED

\begin{claim}\label{cla-sum}
 $\sum\limits_{u\in V(H)}(d_H(u)-1)\frac{x_u}{x_{u^*}}<e(H)+\frac{2\lambda-3}{2\lambda(\lambda-3)}e(H,W_0)+\frac{2\lambda-3}{(2\lambda-4)(\lambda-3)}e(W)$ for any $H\in \mathcal{H}_+$.
 \end{claim}
\noindent {Proof of Claim \ref{cla-sum}}. Recall that for any $H\in \mathcal{H}_+$, we have $H$ is isomorphic to $K_4$ or $K_4-e$. Hence, we will prove the claim by considering the following two cases.

\textbf{Case 1}. $H=K_4$ with the vertex set $V(H)=\{v_1,v_2,v_3,v_4\}$.

In this case, $e(H)=6$, and  $N_{W}(v_i)\cap N_{W}(v_j)=\varnothing$ for any $1\leq i\neq j\leq 4$ since $G^*$ is $F_6$-free. Moreover, $\sum_{i=1}^{4}\sum_{w\in N_{W\setminus W_0}(v_i)}x_w\leq \sum_{w\in W\setminus W_0}x_w\leq 2e(W)x_{u^*}$.
From the eigenequations of vertices in $V(H)$, we obtain
\begin{align*}
(\lambda-3)\sum\limits_{i=1}^{4}x_{v_i}=4x_{u^*}+\sum\limits_{i=1}^{4}\sum\limits_{w\in N_W(v_i)}x_w.
\end{align*}
Hence, by Claim \ref{cla-w0},
\begin{align*}
\sum\limits_{i=1}^{4}x_{v_i}&=\frac{4}{\lambda-3}x_{u^*}+\sum\limits_{i=1}^{4}\sum\limits_{w\in N_W(v_i)}\frac{x_w}{\lambda-3}\\
&=\frac{4}{\lambda-3}x_{u^*}+\sum\limits_{i=1}^{4}\sum\limits_{w\in N_{W_0}(v_i)}\frac{x_w}{\lambda-3}+\sum\limits_{i=1}^{4}\sum\limits_{w\in N_{W\setminus W_0}(v_i)}\frac{x_w}{\lambda-3}\\
&\leq \frac{4}{\lambda-3}x_{u^*}+\frac{2\lambda-3}{2\lambda(\lambda-3)}e(H,W_0)x_{u^*}+\frac{2e(W)}{\lambda-3}x_{u^*}.
\end{align*}
Besides, by (\ref{eqb}), Lemma \ref{le-etah} and Claim \ref{cla-h1}, we have $e(W)< 2|\mathcal{H}_{+}|+\frac{3}{2}<\frac{4\lambda^2-4\lambda+6}{5\lambda-11}+\frac{3}{2}$. Moreover, we get $4+\frac{4}{\lambda-3}+\frac{2\lambda-5}{(2\lambda-4)(\lambda-3)}e(W)<6$ as  $\lambda>55$. Therefore,
\begin{align*}
&\sum_{u\in V(H)}(d_{H}(u)-1)\frac{x_u}{x_{u^*}}=2\sum\limits_{i=1}^{4}\frac{x_{v_i}}{x_{u^*}}\\
&\leq 4+\frac{4}{\lambda-3}+\frac{2\lambda-3}{2\lambda(\lambda-3)}e(H,W_0)+\frac{2e(W)}{\lambda-3}\\
&=4+\frac{4}{\lambda-3}+\frac{2\lambda-5}{(2\lambda-4)(\lambda-3)}e(W)+\frac{2\lambda-3}{2\lambda(\lambda-3)}e(H,W_0)+\frac{2\lambda-3}{(2\lambda-4)(\lambda-3)}e(W)\\
&<6+\frac{2\lambda-3}{2\lambda(\lambda-3)}e(H,W_0)+\frac{2\lambda-3}{(2\lambda-4)(\lambda-3)}e(W).
\end{align*}

\textbf{Case 2}. $H=K_4-e$ with the vertex set $V(H)=\{v_1,v_2,v_3,v_4\}$, where $d_{H}(v_2)=d_{H}(v_3)=3$.

In this case, $e(H)=5$ and $N_{W}(v_i)\cap N_{W}(v_j)=\varnothing$ for any $1\leq i\neq j\leq 3$ since $G^*$ is $F_6$-free. The eigenequations of $v_1$, $v_2$ and $v_3$ are as follows.
\begin{align*}
\left\{
\begin{array}{ll}
\lambda x_{v_1} &= x_{v_2}+x_{v_3}+x_{u^*}+\sum\limits_{w\in N_W(v_1)}x_w;\\
\lambda x_{v_2} &= x_{v_1}+x_{v_3}+x_{v_4}+x_{u^*}+\sum\limits_{w\in N_W(v_2)}x_w;\\
\lambda x_{v_3} &= x_{v_1}+x_{v_2}+x_{v_4}+x_{u^*}+\sum\limits_{w\in N_W(v_3)}x_w.
\end{array}
\right.
\end{align*}
Adding the above three equations together, we get
\begin{align*}
(\lambda-2)\sum\limits_{i=1}^{3}x_{v_i}\leq 5 x_{u^*}+\sum\limits_{i=1}^{3}\sum\limits_{w\in N_W(v_i)}x_w.
\end{align*}
 Similarly to the proof of Case 1, we can obtain
 \begin{align*}
\sum\limits_{i=1}^{3}x_{v_i}\leq \frac{5}{\lambda-2}x_{u^*}+\frac{2\lambda-3}{2\lambda(\lambda-2)}e(H,W_0)x_{u^*}+\frac{2e(W)}{\lambda-2}x_{u^*}.
\end{align*}
Furthermore,
\begin{align*}
\sum_{u\in V(H)}(d_{H}(u)-1)\frac{x_u}{x_{u^*}}&=\frac{2x_{v_2}+2x_{v_3}+x_{v_1}+x_{v_4}}{x_{u^*}}\\
&\leq 3+\frac{5}{\lambda-2}+\frac{2\lambda-3}{2\lambda(\lambda-2)}e(H,W_0)+\frac{2e(W)}{\lambda-2}\\
&=3+\frac{5}{\lambda-2}+\frac{1}{\lambda-2}e(W)+\frac{2\lambda-3}{2\lambda(\lambda-2)}e(H,W_0)+\frac{1}{\lambda-2}e(W)\\
&<5+\frac{2\lambda-3}{2\lambda(\lambda-3)}e(H,W_0)+\frac{2\lambda-3}{(2\lambda-4)(\lambda-3)}e(W),
\end{align*}
where the last inequality derives from $3+\frac{5}{\lambda-2}+\frac{1}{\lambda-2}e(W)<5$ as $\lambda>55$ and $e(W)< \frac{4\lambda^2-4\lambda+6}{5\lambda-11}+\frac{3}{2}$.
Thus Claim \ref{cla-sum} holds.\QED

Now we come back to prove Lemma \ref{le-etah0}. For any $H\in \mathcal{H}\setminus \mathcal{H}_{+}$, $\eta(H)\leq 0$ implies that $\sum_{u\in V(H)}(d_H(u)-1)\frac{x_{u}}{x_{u^*}}\leq e(H)$.  In the following inequalities, we note that  the first inequality comes from Claims \ref{cla-w0} and \ref{cla-sum}, the second inequality results from Claim \ref{cla-h1}, the third inequality derives from $\lambda>55$ and the last inequality deduces from $e(W)\geq 4$.
\begin{align*}
&\sum\limits_{u\in U_{+}}(d_{U_{+}}(u)-1)\frac{x_{u}}{x_{u^*}}+\sum\limits_{w\in W}d_U(w)\frac{x_w}{x_{u^*}}\\
&=\sum\limits_{H\in \mathcal{H}_{+}}\sum\limits_{u\in V(H)}(d_H(u)-1)\frac{x_{u}}{x_{u^*}}+\sum\limits_{H\in \mathcal{H}\setminus\mathcal{H}_{+}}\sum\limits_{u\in V(H)}(d_H(u)-1)\frac{x_{u}}{x_{u^*}}\\
&+\sum\limits_{w\in W_{0}}d_U(w)\frac{x_w}{x_{u^*}}+\sum\limits_{w\in W\setminus W_{0}}d_U(w)\frac{x_w}{x_{u^*}}\\
&<\sum\limits_{H\in \mathcal{H}_{+}}\left(e(H)+\frac{2\lambda-3}{2\lambda(\lambda-3)}e(H,W_0)+\frac{2\lambda-3}{(2\lambda-4)(\lambda-3)}e(W)\right)+\sum\limits_{H\in \mathcal{H}\setminus\mathcal{H}_{+}}e(H)\\
&+\frac{2\lambda-3}{2\lambda}e(U,W_0)
+e(U,W\setminus W_0)\\
&<e(U_{+})+\left(\frac{2\lambda-3}{2\lambda(\lambda-3)}+\frac{2\lambda-3}{2\lambda}\right)e(U,W_0)+e(U,W\setminus W_0)\\
&+\frac{(2\lambda-3)(2\lambda^2-2\lambda +3)}{(2\lambda-4)(\lambda-3)(5\lambda-11)}e(W)\\
&<e(U_{+})+e(U,W)+\frac{5}{8}e(W)\\
&\leq e(U_{+})+e(U,W)+e(W)-\frac{3}{2}.
\end{align*}
Hence, we get a contradiction to (\ref{eq-etb}). The proof is complete.\QED

Recall that  $\gamma (w)=d_U(w)(1-\frac{x_{w}}{x_{u^*}})$ for any $w\in W$. By Lemma \ref{le-etah0} and (\ref{eqb}), we have $e(W)\leq 1$, $\eta(H)> -\frac{3}{2}$ holds for any $H\in \mathcal{H}$ and $\gamma(w)<\frac{3}{2}$ holds for any $w\in W$.

\begin{lemma}\label{le-k4}
$H\neq K_4$ for any $H\in \mathcal{H}$.
\end{lemma}
\noindent\textbf{Proof.} Suppose there exists a component $H'=K_4\in \mathcal{H}$ with $V(H')=\{u_1,u_2,u_3,u_4\}$. Without loss of generality, let $x_{u_1}\geq x_{u_2}\geq x_{u_3}\geq x_{u_4}$. Note that $N_W(u_i)\cap N_W(u_j)=\varnothing$ holds for any $1\leq i\neq j\leq 4$ since $G^*$ is $F_6$-free.

If $W_{H'}=\varnothing$, then $x_{u_1}=x_{u_2}=x_{u_3}=x_{u_4}=\frac{x_{u^*}}{\lambda-3}$. Hence, as $\lambda>55$, $\eta(H')=2 \sum_{i=1}^{4}\frac{x_{u_i}}{x_{u^*}}-6=\frac{8}{\lambda-3}-6\leq-\frac{3}{2}$. It is a contradiction.

If $W_{H'}\neq \varnothing$, then we claim that $N_W(u_2)=N_W(u_3)=N_W(u_4)=\varnothing$. Otherwise we can construct a new $F_6$-free graph with $m$ edges and larger spectral radius as $e(W)\leq 1$. It follows that $x_{u_2}=x_{u_3}=x_{u_4}\leq \frac{2x_{u^*}}{\lambda-2}$.
Hence, as $\lambda>55$, $\eta(H')=2x_{u_1}+2 \sum_{i=2}^{4}\frac{x_{u_i}}{x_{u^*}}-6\leq\frac{12}{\lambda-2}-4<-\frac{3}{2}$. It is a contradiction.\QED

\begin{lemma}\label{le-k4e}
$H\neq K_4-e$ for any $H\in \mathcal{H}$.
\end{lemma}
\noindent\textbf{Proof.} Suppose there exists a component $H'=K_4-e\in \mathcal{H}$ with vertex set $V(H')=\{u_1,u_2,u_3,u_4\}$, where $d_{H'}(u_2)=d_{H'}(u_3)=3$. Without loss of generality, let $ x_{u_2}\geq x_{u_3}$. Note that $N_W(u_i)\cap N_W(u_j)=\varnothing$ holds for any $1\leq i\neq j\leq 3$ since $G^*$ is $F_6$-free.

If $N_W(u_2)=N_W(u_3)=\varnothing$, then $x_{u_2}=x_{u_3}=\frac{3x_{u^*}}{\lambda-1}$. Hence, as $\lambda>55$, $\eta(H')= \frac{2x_{u_2}+2x_{u_3}+x_{u_1}+x_{u_4}}{x_{u^*}}-5\leq\frac{12}{\lambda-1}-3<-\frac{3}{2}$ . It is a contradiction.

If $N_W(u_2)\neq \varnothing $ or $N_W(u_3)\neq \varnothing$, then we claim that $N_W(u_3)=\varnothing$. Otherwise we can construct a new $F_6$-free graph with $m$ edges and larger spectral radius as $e(W)\leq 1$. Therefore, $x_{u_3}\leq \frac{4x_{u^*}}{\lambda}$. Besides, $W_0\neq \varnothing$. Otherwise, since $e(W)\leq 1$, we obtain $|W|\leq 2$. It follows that $x_{u_2}\leq\frac {6x_{u^*}}{\lambda}$. Furthermore, $\eta(H')\leq -\frac{3}{2}$ as $\lambda>55$, a contradiction. \textcolor{black}{Similarly, we can get $|N_{W_0}(u_2)|\geq 1$.}

\begin{claim}\label{cla-w01}
There exists an vertex $w_0\in W_0$ such that $x_{w_0}>\frac{\lambda-1}{\lambda}x_{u^*}$.
\end{claim}

\noindent{Proof of Claim \ref{cla-w01}}. Suppose to the contrary that $x_{w}\leq \frac{\lambda-1}{\lambda}x_{u^*}$   holds for any $w\in W_0$.  Since $G^*$ is $F_6$-free, we have $N_{W}(u_1)\cap N_{W}(u_2)=\varnothing$. Besides,
\begin{align*}
\left\{
\begin{array}{ll}
\lambda x_{u_1}&= x_{u_2}+x_{u_3}+x_{u^*}+\sum\limits_{w\in N_W(u_1)}x_w;\\
\lambda x_{u_2}&= x_{u_1}+x_{u_3}+x_{u_4}+x_{u^*}+\sum\limits_{w\in N_W(u_2)}x_w.
\end{array}
\right.
\end{align*}
Hence,
\begin{align}
x_{u_1}+x_{u_2}&=\frac{2x_{u^*}+2x_{u_3}+x_{u_4}+\sum\limits_{w\in N_W(u_1)\cup N_W(u_2)}x_w}{\lambda-1}\nonumber\\
&=\frac{2x_{u^*}+2x_{u_3}+x_{u_4}}{\lambda-1}+\sum\limits_{w\in N_{W_0}(u_1)\cup N_{W_0}(u_2)}\frac{x_w}{\lambda-1}+\sum\limits_{w\in N_{W\setminus W_0}(u_1)\cup N_{W\setminus W_0}(u_2)}\frac{x_w}{\lambda-1}\nonumber\\
&\leq \left(\frac{3}{\lambda-1}+\frac{8}{\lambda(\lambda-1)}\right)x_{u^*}+\frac{e(H',W_0)}{\lambda}x_{u^*}+\frac{2e(W)}{\lambda-1}x_{u^*}.\label{eqd}
\end{align}
Then,  we can get
\begin{align}
&\sum\limits_{u\in U_{+}}(d_U(u)-1)\frac{x_{u}}{x_{u^*}}+\sum\limits_{w\in W}d_U(w)\frac{x_w}{x_{u^*}} \nonumber\\
&=\sum\limits_{u\in V(H')}(d_{H'}(u)-1)\frac{x_{u}}{x_{u^*}}+\sum\limits_{H\in \mathcal{H}\setminus H' }\sum\limits_{u\in V(H)}(d_H(u)-1)\frac{x_{u}}{x_{u^*}}\\
&+\sum\limits_{w\in W_{0}}d_U(w)\frac{x_w}{x_{u^*}}+\sum\limits_{w\in W\setminus W_{0}}d_U(w)\frac{x_w}{x_{u^*}} \nonumber \\
&\leq2+\frac{8}{\lambda}+\frac{3}{\lambda-1}+\frac{8}{\lambda(\lambda-1)}+\frac{e(H',W_0)}{\lambda}+\frac{2e(W)}{\lambda-1}\\
&+\sum\limits_{H\in \mathcal{H}\setminus H'}e(H)+\frac{\lambda-1}{\lambda}e(U,W_0)+e(U,W\setminus W_0)\nonumber\\
&< e(U_{+})+e(U,W)+e(W)-\frac{3}{2}, \label{eq-con1}
\end{align}
where the second inequality results from Lemma \ref{le-etah0} and (\ref{eqd}), and the last inequality holds since $e(H')=5$ and $\lambda>55$. We can see that (\ref{eq-con1})  contradicts to (\ref{eqc}). This completes the proof. \QED

We proceed by distinguishing the following two cases.

\textbf{Case 1}. $x_{u_2}+x_{u_3}\geq x_{u^*}$.

In this case, we claim that $x_{w}\leq \frac{\lambda-1}{\lambda}x_{u^*}$ holds for any $w\in W_0$. Indeed, for  a vertex $w\in W_0$, if $w$ is adjacent to $u_2$, then $w$ is not adjacent to $u_1$, $u_3$ and $u_4$. Moreover, $\lambda x_w\leq \lambda x_{u^*}-x_{u_1}-x_{u_3}-x_{u_4}$. By Lemma \ref{le-etah0} and (\ref{eqb}), we know that $\eta(H')=\frac{2x_{u_2}+2x_{u_3}+x_{u_1}+x_{u_4}}{x_{u^*}}-5>-\frac{3}{2}$. Recall that $x_{u_3}\leq \frac{4}{\lambda}x_{u^*}$, we have $x_{u_3}+x_{u_1}+x_{u_4}>x_{u^*}$ as $\lambda>55$. Thus, $\lambda x_w<\lambda x_{u^*}-x_{u^*}$, which implies that $x_{w}\leq \frac{\lambda-1}{\lambda}x_{u^*}$. If $w$ is not adjacent to $u_2$, then $\lambda x_w\leq \lambda x_{u^*}-x_{u_2}-x_{u_3}$. It follows that $x_{w}\leq \frac{\lambda-1}{\lambda}x_{u^*}$. Consequently, $x_w\leq \frac{\lambda-1}{\lambda}x_{u^*}$ holds for any $w\in W_0$, which contradicts to Claim \ref{cla-w01}.

\textbf{Case 2}. $x_{u_2}+x_{u_3}< x_{u^*}$.

 In this case, since $\eta(H')>-\frac{3}{2}$, we have that $x_{u_1}+x_{u_4}>\frac{3}{2}x_{u^*}$.
Recall that  $|N_{W_0}(u_2)|\geq 1$. Let $w_1\in N_{W_0}(u_2)$. By Lemma \ref{le-zhai}, we have $d_U(w_1)>2$. It follows that $x_{w_1}>\frac{1}{4}x_{u^*}$ from $\gamma(w_1)<\frac{3}{2}$.
Then we can construct a graph $G'=G^*-\{u_2u_3,u_2w_1\}+\{w_1u_1,w_1u_4\}$. Clearly, $G'$ is $F_6$-free and $\lambda(G')>\lambda$. Indeed,
\begin{align*}
\lambda(G')-\lambda&\geq \mathbf{x}^{\mathrm{T}}A(G')\mathbf{x}-\mathbf{x}^{\mathrm{T}}A(G)\mathbf{x}\\
&=2(x_{w_1}x_{u_1}+x_{w_1}x_{u_4})-2(x_{u_3}x_{u_2}+x_{w_1}x_{u_2})\\
&>x_{w_1}x_{u^*}-2x_{u_2}x_{u_3}\\
&>(\frac{1}{4}-\frac{2}{\lambda})x_{u^*}^2\\
&>0,
\end{align*}
where the last inequality comes from $\lambda>55$.
It is a contradiction to the choice of $G^*$.

To sum up, Lemma \ref{le-k4e} holds.\QED

Recall that  $\gamma (w)=d_U(w)(1-\frac{x_{w}}{x_{u^*}})$ for any $w\in W$ and $\gamma(W')=\sum_{w\in W'}\gamma(w)$ for any non-empty subset $W'\subseteq W$.

\begin{lemma}\label{le-ew0}
$e(W)=0$.
\end{lemma}
\noindent\textbf{Proof.} Suppose to the contrary, $e(W)=1$. By Lemma \ref{le-zhai}, we know that $d_U(w)\geq 1$ for any $w\in W$.  And by (\ref{eqb}), we will see that $\eta(H)>-\frac{1}{2}$ holds for any $H\in \mathcal{H}$, $\sum_{u\in U_0}\frac{x_u}{x_{u^*}}<\frac{1}{2}$ and $\gamma(W)<\frac{1}{2}$. From $\gamma(W)<\frac{1}{2}$, we can get that $\gamma(W')<\frac{1}{2}$ for any non-empty subset $W'\subseteq W$. Hence, for any $w\in W$, $x_w>\frac{1}{2}x_{u^*}$, which yields that $d_{U}(w)> \frac{\lambda}{2}$ for any $w\in W_0$.
Since $\eta(H)>-\frac{1}{2}$ holds for any $H\in \mathcal{H}$, we know that $x_u>\frac{1}{2}x_{u^*}$ for any $u\in V(H)$ with $d_{H}(u)\geq 2$. Moreover, by Lemmas \ref{le-etah}, \ref{le-k4} and \ref{le-k4e}, every element of $\mathcal{H}$ is isomorphic to a $C_4$ or a $K_{1,r}+e$ for some positive integer $r$.  Therefore, we distinguish the following two cases to lead a contradiction,  respectively.

\textbf{Case 1}. There is a component $H'=C_4\in \mathcal{H}$.

Since for any $u\in V(H')$, $x_u>\frac{1}{2}x_{u^*}$. We have $d_W(u)\geq 5$ as $\lambda>55$. Furthermore, $|W_0|\geq 3$. However, there are not two vertices in $W_0$ with degrees at least 3 due to that $G^*$ is $F_6$-free. Hence, we know that  there exists two vertices $w_1$ and $w_2$ in $W_0$ such that $d_{H'}(w_1)\leq 2$ and $d_{H'}(w_2)\leq 2$. It follows that  $x_{w_1}, x_{w_2}< \frac{\lambda-1}{\lambda}x_{u^*}$ since $x_u>\frac{1}{2}x_{u^*}$ for any $u\in V(H')$.
Conclusively, let $W'=\{w_1,w_2\}$, we get
\begin{align*}
\gamma(W')&=\sum\limits_{i=1}^{2}d_U(w_i)(1-\frac{x_{w_i}}{x_{u^*}})>2\cdot \frac{\lambda}{2}(1-\frac{\lambda-1}{\lambda})
=\frac{1}{2},
\end{align*}
which is a contradiction.

 \textbf{Case 2}. There is a component $H'=K_{1,r}+e\in \mathcal{H}$ for some positive integer $r\geq 2$.

 Assume the set of vertices in $H'$ which form the triangle is  $T=\{u_1,u_2,u_3\}$, and the center vertex is $u_1$. Similar to the above case, since $x_{u_2}, x_{u_3}>\frac{1}{2}x_{u^*}$, we know that $|W_0|\geq 4$ as $\lambda>55$.
We will show that there is no vertex  $w\in W_0$ such that $T\subseteq N(w)$. Otherwise, suppose there is a vertex $w_0\in W_0$ and $T\subseteq N(w_0)$, then we know that $|N_T(w)|\leq 1$ for any $w\in W_0\setminus \{w_0\}$.  Furthermore, $x_w< \frac{\lambda-1}{\lambda}x_{u^*}$ as $x_{u_i}>\frac{1}{2}x_{u^*}$ for any $u_i\in T$, which implies $\gamma(w)>\frac{1}{2}$ as $d_{U}(w)>\frac{\lambda}{2}$. It is a contradiction.
Hence, for any $w\in W_0$, $|N_T(w)|\leq 2$, which yields $x_w< \frac{2\lambda-1}{2\lambda}x_{u^*}$. Now we could choose any subset $W'\subseteq W_0$ with $|W'|=4$, and see that $\gamma(W')>\frac{1}{2}$, which implies a contradiction. Indeed, since $x_{u_i}>\frac{1}{2}x_{u^*}$ for any $u_i\in T$,
\begin{align*}
\gamma(W')=\sum\limits_{w\in W'}d_U(w)(1-\frac{x_{w}}{x_{u^*}})>4\cdot \frac{\lambda}{2}(1-\frac{2\lambda-1}{2\lambda})=\frac{1}{2}.
\end{align*}

In summary, $\mathcal{H}=\varnothing$, which contradicts to Lemma \ref{le-eU}.\QED

By Lemma \ref{le-ew0}, $W=W_0$. Combining it with  Lemma \ref{le-zhai},  we have $d_{U}(w)\geq 2$ for any $w\in W$. Recall that by (\ref{eqb}),  it is obvious that $\eta(H)>-\frac{3}{2}$ holds for any $H\in \mathcal{H}$, $\sum_{u\in U_0}\frac{x_u}{x_{u^*}}<\frac{3}{2}$ and $\gamma(W)<\frac{3}{2}$. From $\gamma(W)<\frac{3}{2}$, we can get that $\gamma(W')<\frac{3}{2}$ for any non-empty subset $W'\subseteq W$.

\begin{lemma}\label{le-wd}
There is at most one vertex in $W$ with degree less than $\frac{\lambda}{2}$.
\end{lemma}
\noindent\textbf{Proof}. Otherwise, suppose that there exist $w$, $w'\in W$ such that $d(w)$, $d(w')< \frac{\lambda}{2}$. Hence, $x_{w}$, $x_{w'}<\frac{1}{2}x_{u^*}$.  It is easy to check that $\gamma(\{w,w'\})>\frac{3}{2}$ since $d_U(w), d_U(w')\geq 2$, a contradiction.\QED

\begin{lemma}\label{le-c4}
 $H\neq C_4$ for any $H\in \mathcal{H}$.
 \end{lemma}
\noindent \textbf{Proof.} Suppose there is a component $H'=C_4\in \mathcal{H}$ with $V(H')=\{u_1,u_2,u_3,u_4\}$.  It should be noted that $x_{u_i}+x_{u_j}>\frac{1}{2}x_{u^*}$ for any $1\leq i\neq j\leq 4$. Indeed, if not, we will have $\eta(H')\leq -\frac{3}{2}$, a contradiction. Furthermore, there is at least one vertex $u\in V(H')$ satisfying $x_u>\frac{1}{4}x_{u^*}$, which implies that $|W|\geq 8$ as $\lambda>55$ as the previous proof. Recall that there are not two vertices in $W$ with degrees at least 3 due to $G^*$ is $F_6$-free.  Hence, combining it with  Lemma \ref{le-wd}, we could choose $W'\subseteq W$ with $|W'|= 6$ such that $d_{H'}(w)\leq 2$  and $d_{U}(w)\geq \frac{\lambda}{2}$ holds for any $w\in W'$. It follows that $x_w\leq \frac{2\lambda-1}{2\lambda}x_{u^*}$ for any $w\in W'$. Then we obtain
\begin{align*}
\gamma(W')=\sum\limits_{w\in W'}d_U(w)(1-\frac{x_{w}}{x_{u^*}})>6\cdot \frac{\lambda}{2}(1-\frac{2\lambda-1}{2\lambda})=\frac{3}{2},
\end{align*}
which is a contradiction.\QED

\begin{lemma}\label{le-k1re}
For any $H\in \mathcal{H}$, $H\neq K_{1,r}+e$, where $r\geq 2$ is a positive integer.
\end{lemma}
\noindent \textbf{Proof.} It is obvious that if there is a component of $\mathcal{H}$ isomorphic to star or double star, then by Lemma \ref{le-etah}, we can get a contradiction by a similar proof to Case 2 in Lemma \ref{le-ew0}.  Together with Lemmas \ref{le-k4}, \ref{le-k4e} and \ref{le-c4}, we know that every component of $\mathcal{H}$ is isomorphic to  $K_{1,r}+e$ for some positive integer $r$.

\begin{claim}\label{cl-h11}
$|\mathcal{H}|=1$.
\end{claim}
\noindent {Proof of Claim \ref{cl-h11}.} Suppose to the contrary that $|\mathcal{H}|\geq 2$. Assume that there are two components $H_1, H_2\in \mathcal{H}$ with $V_1\subseteq V(H_1)$ and $V_2\subseteq V(H_2)$, where the vertices of $V_1$ and $V_2$ form the triangles in $H_1$ and $H_2$, respectively. For $i=1,2$, we could easily get that  there is at most one vertex $w_i\in W$ such that $V_i\subseteq N_{H_i}(w_i)$ as $G^*$ is $F_6$-free.  On the other hand, similarly to the previous proof, we could prove that $|W|\geq 9$ as $\lambda>55$.  Hence, we could find a subset $W\subseteq W'$ with $|W'|=6$ such that $d_{V_1}(w)\leq 2$, $d_{V_2}(w)\leq 2$ and $d_{U}(w)>\frac{\lambda}{2}$ hold for any $w\in W'$.
Note that for any two  distinct vertices $u,v\in V_1\cup V_2$, we have $x_u+x_v>\frac{1}{2}x_{u^*}$ since $\eta(H_1)+\eta(H_2)>-\frac{3}{2}$.
Consequently,
\begin{align*}
\gamma(W')=\sum\limits_{w\in W'}d_U(w)(1-\frac{x_{w}}{x_{u^*}})>6\cdot \frac{\lambda}{2}(1-\frac{2\lambda-1}{2\lambda})=\frac{3}{2},
\end{align*}
which is a contradiction.\QED

Suppose the unique component of $\mathcal{H}$  is $H'=K_{1,r}+e$ with $V=\{u_1,u_2,u_3\}\subseteq V(H')$, where $u_1$ is the center vertex and $d_{H'}(u_2)=d_{H'}(u_3)=2$.
 We claim that $r\geq 3$. Otherwise, $\lambda x_{u^*}\leq 3x_{u^*}+\sum_{u\in U_0}x_u$, which yields that $\lambda<\frac{9}{2}$ from  $\sum_{u\in U_0}x_u< \frac{3}{2}x_{u^*}$. It is a contradiction as $\lambda>55$.

 If there exists a vertex $w_0\in  N_W(u_1)$, then $w_0$ is not adjacent to any vertex of $V(H')\setminus\{u_1\}$  since $G^*$ is $F_6$-free and $r\geq 3$.
 It follows that $x_{w_0}\leq \frac{1}{\lambda} (x_{u^*}+\sum_{u\in U_0}x_u)<\frac{5 x_{u*}}{2\lambda}$, which implies $\gamma(w_0)=d_U(w_0)(1-\frac{5 x_{u*}}{2\lambda})>\frac{3}{2}$, a contradiction  as $d_{U}(w_0)\geq 2$ and $\lambda>55$.
Hence, $N_W(u_1)= \varnothing$.
If $r=3$, then $x_{u_1}\leq \frac{4}{\lambda}x_{u^*}$. It is a contradiction to $\lambda>55$. Thus, $r\geq 4$. Furthermore, $x_{u_1}>\frac{1}{2}x_{u^*}$ as $r\geq 4$ and $\eta(H')>-\frac{3}{2}$.
For any $w\in W$, $x_w< \frac{\lambda x_{u^*}-x_{u_1}}{\lambda}<\frac{2\lambda-1}{2\lambda}x_{u^*}$.
Similarly to the previous proof, we could prove that $|W|\geq 7$ as $\lambda>55$.  Hence, we could find a subset $W\subseteq W'$ with $|W'|=6$ such that $d_U(w)>\frac{\lambda}{2}$ holds for any $w\in W'$. Then we could see that
 \begin{align*}
\gamma(W')=\sum\limits_{w\in W'}d_U(w)(1-\frac{x_{w}}{x_{u^*}})>6\cdot \frac{\lambda}{2}(1-\frac{2\lambda-1}{2\lambda})=\frac{3}{2},
\end{align*}
which is a contradiction.\QED

By Lemmas \ref{le-k4}, \ref{le-k4e}, \ref{le-c4} and \ref{le-k1re}, we get that every component of $\mathcal{H}$ is a star  or a double star. Combining it with Lemma \ref{le-etah} and (\ref{eqb}), we have that $|\mathcal{H}|=1$.  Hence, let $\mathcal{H}=\{H\}$. Furthermore, we will see that $\eta(H)>-\frac{1}{2}$, $\sum_{u\in U_0}\frac{x_u}{x_{u^*}}<\frac{1}{2}$ and $\gamma(W)<\frac{1}{2}$. From $\gamma(W)<\frac{1}{2}$, we can get that $\gamma(W')<\frac{1}{2}$ for any non-empty subset $W'\subseteq W$.

\begin{lemma}
$H$ is a star.
\end{lemma}
\noindent\textbf{Proof.}
It remains to prove that the unique component is not a double star. Suppose to the contrary that the unique component is a double star $H$ with center vertices $u_0$ and $v_0$. Let $V_1=\{u_1,\ldots, u_a\}$ and $V_2=\{v_1,\ldots, v_b\}$ be the leaf vertices of $H$ adjacent to $u_0$ and $v_0$, respectively, where $a, b\geq 1$ are positive integers. Since $\eta(H)>-\frac{1}{2}$, we see that $x_{u_0},x_{v_0}>\frac{1}{2}x_{u^*}$. Similarly, from $\gamma(W')<\frac{1}{2}$ for any $W'\subseteq W$, we could get
$x_w>\frac{1}{2}x_{u^*}$ and  $d_U(w)>\frac{\lambda}{2}$ for  any $w\in W$.

Next  we will prove this lemma by considering the following two cases.

\textbf{Case 1.} $W_H=\varnothing$.

 In this case, the eigenequations of $u^*$, $u_0$ and $v_0$ are  as follows.
\begin{align*}
\left\{
 \begin{array}{ll}
 \lambda x_{u^*}&=\sum\limits_{i=0}^{a}x_{u_i}+\sum\limits_{i=0}^{b}x_{v_i}+\sum\limits_{u\in U_0}x_u;\\
 \lambda x_{u_0}&=x_{u^*}+x_{v_0}+\sum\limits_{i=1}^{a}x_{u_i};\\
 \lambda x_{v_0}&=x_{u^*}+x_{u_0}+\sum\limits_{i=1}^{b}x_{v_i}.
 \end{array}
\right.
\end{align*}
In fact, we have  $a, b\geq 2$ since $\lambda>55$ and $x_{u_0},x_{v_0}>\frac{1}{2}x_{u^*}$.  Furthermore, by using $\eta(H)>-\frac{1}{2}$ and $a, b\geq 2$, we could immediately get $x_{u_0},x_{v_0}>\frac{3}{4}x_{u^*}$. Thus,
\begin{align*}
\frac{3}{2}\lambda x_{u^*}<\lambda (x_{u_0}+x_{v_0})=\lambda x_{u^*}+2x_{u^*}-\sum\limits_{u\in U_0}x_u,
\end{align*}
which is a contradiction as $\lambda>55$.

\textbf{Case 2.} $W_H\neq \varnothing$.

Before discussing Case 2, we need to show that $|N(w)\cap \{u_0,v_0\}|\geq 1$ holds for any $w\in W$. In fact, if there exists a vertex $w\in W$ such that $N(w)\cap \{u_0,v_0\}=\varnothing$, then $\lambda x_{w}\leq \lambda x_{u^*}-x_{u_0}-x_{v_0}<(\lambda-1)x_{u^*}$ as $x_{u_0},x_{v_0}>\frac{1}{2}x_{u^*}$. Then it is  easy to verify that $\gamma(w)>\frac{1}{2}$ as $d_{U}(w)>\frac{\lambda}{2}$, a contradiction.

Now, we proceed by distinguishing two subcases.

\textbf{Subcase 2.1.} $a,b\geq 2$.

In this case, we could easily get that $x_{u_0},x_{v_0}>\frac{3}{4}x_{u^*}$ since $\eta(H)>-\frac{1}{2}$.
We firstly claim that $|N(w)\cap \{u_0,v_0\}|=1$ holds for any $w\in W$. Otherwise, suppose that there exists a vertex $w\in W$ such that $\{u_0,v_0\}\subseteq N(w)$. Then we know that $|N_{V_1}(w)|\leq 1$ and  $|N_{V_2}(w)|\leq 1$. Hence, $\frac{\lambda}{2}x_{u^*}<\lambda x_w\leq 4x_{u^*}+\sum_{u\in  U_0}x_u<\frac{9}{2}x_{u^*}$, which implies that $\lambda<9$, a contradiction as $\lambda>55$. Hence,  for any $w\in W$, $|N(w)\cap \{u_0,v_0\}|=1$ holds,  and it follows that  $x_{w}\leq \frac{4\lambda-3}{4\lambda}x_{u^*}$. If  $|W_H|\geq 2$, then  we could select any subset $W'\subseteq W$ with $|W'|=2$. It is easy to check that
$\gamma(W')>\frac{1}{2}$, a contradiction. Therefore, $|W_H|=1$.  Assume that $W_H=\{w\}$ and $w$ is adjacent to $u_0$. Hence, $|N_{V_1}(w)|\leq 1$. Without loss of generality, let $u_i\notin N(w)$ for $i=2,\ldots,a$. Hence,
\begin{align*}
\left\{
 \begin{array}{ll}
 \lambda x_{u^*}&=\sum\limits_{i=0}^{a}x_{u_i}+\sum\limits_{i=0}^{b}x_{v_i}+\sum\limits_{u\in U_0}x_u;\\
 \lambda x_{u_0}&=x_{u^*}+x_{v_0}+\sum\limits_{i=1}^{a}x_{u_i}+x_w.
 \end{array}
\right.
\end{align*}
Moreover, since $x_{u_0}>\frac{3}{4}x_{u^*}$ and $x_{w}\leq \frac{4\lambda-3}{4\lambda}x_{u^*}$,
\begin{align*}
\lambda x_{w}\leq \lambda x_{u^*}-\sum\limits_{i=2}^{a}x_{u_i}= \lambda x_{u^*}-\lambda x_{u_0}+x_{v_0}+x_{u^*}+x_{w}\leq \frac{\lambda^2+12\lambda-3}{4\lambda}x_{u^*},
\end{align*}
which implies that $x_{w}\leq \frac{\lambda^2+12\lambda-3}{4\lambda^2}x_{u^*}$. Thus, $\gamma(w)>\frac{1}{2}$ as $\lambda>55$, a contradiction.

\textbf{Subcase 2.2.} $a=1$ or $b=1$.

 If  $a=1$ and $b=1$, then we know that $\lambda x_w\leq 4 x_{u^*}+\sum_{u\in U_0}x_u<\frac{9}{2}x_{u^*}$ for any $w\in W$, which is a contradiction to $\lambda>55$. Hence, without loss of generality, we can assume that $a\geq 2$ and $b=1$. We will prove that $|N(w)\cap \{u_0,v_0\}|=1$ for any $w\in W$. Otherwise, suppose that $w\in W$ and $\{u_0,v_0\}\subseteq N(w)$. Then  $|N_{V_1}(w)|\leq 1$. It follows that $\lambda x_w\leq 4 x_{u^*}+\sum_{u\in U_0}x_u<\frac{9}{2}x_{u^*}$, which is a contradiction to $\lambda>55$. Therefore, $|N(w)\cap \{u_0,v_0\}|=1$ for any $w\in W$. Furthermore, $x_{w}<\frac{2\lambda-1}{2\lambda}x_{u^*}$.
 Note that $|W|\geq 2$ since $x_{v_0}\geq \frac{1}{2}x_{u^*}$ and $\lambda>55$. Now we could arbitrarily choose a subset $W'\subseteq W$ with $|W'|=2$ and $\gamma(W')>\frac{1}{2}$, a contradiction. The proof is complete.\QED

 \begin{lemma}\label{le-wemp}
 $G^*=S_{\frac{m+4}{2},2}^{-1}$.
 \end{lemma}
 \noindent\textbf{Proof.} By Lemma \ref{le-Fang}, it is sufficient to prove that $W=\varnothing$. Otherwise, assume that $W\neq \varnothing$.
 Suppose that $G^*[U]=K_{1,a}\cup bK_1$, where $V(bK_1)=\{v_1,\ldots,v_b\}$ and $V(K_{1,a})=\{u_0,u_1,\ldots,u_a\}$ with the center vertex $u_0$.
 Since $\lambda x_{u^*}=\sum_{i=0}^{a}x_{u_i}+\sum_{i=1}^{b}x_{v_i}$ and $\sum_{i=1}^{b}x_{v_i}<\frac{1}{2}x_{u^*}$, we have $a\geq 3$ as $\lambda>55$.
  \begin{claim}\label{cl-wu}
  $w$ is  not adjacent to $u_0$ for any $w\in W$.
  \end{claim}
  \noindent{Proof of Claim \ref{cl-wu}.} Suppose to the contrary there exists a vertex $w\in W$ such that $w$ is adjacent to $u_0$. Then since $G^*$ is $F_6$-free and $a\geq 3$, we obtain $|N(w)\cap (V(K_{1,a})\setminus \{u_0\})|\leq 1$. Hence, $\lambda x_w\leq 2 x_{u^*}+\sum_{i=1}^{b}x_{v_i}<\frac{5}{2}x_{u^*}$, which is a contradiction to $\lambda>55$. \QED

  \begin{claim}\label{cl-wvi}
  $w$ is  not adjacent to $v_i$ for any $w\in W$ and any $i\in [b]$.
  \end{claim}
  \noindent{Proof of  Claim \ref{cl-wvi}.} Suppose to the contrary there exists a vertex $w'\in W$ such that $w'$ is adjacent to $v_i$.  From the eigenequations of $w'$ and $u_0$, we have
  \begin{align*}
\left\{
 \begin{array}{ll}
 \lambda x_{u_0}&=x_{u^*}+\sum\limits_{i=1}^{a}x_{u_i}+\sum\limits_{w\in N_W(u_0)}x_w;\\
 \lambda x_{w'}&\leq \sum\limits_{i=1}^{a}x_{u_i}+\sum_{i=1}^{b}x_{v_i}.
 \end{array}
\right.
\end{align*}
Combining it with $\sum_{i=1}^{b}x_{v_i}<\frac{1}{2}x_{u^*}$, we  get $x_{w'}< x_{u_0}$. Let $G'=G^*-v_jw'+v_ju_0$. It is obvious that $G^*$ is $F_6$-free, and $\lambda(G')>\lambda$ by Lemma \ref{le-zhai}, which is a contradiction.\QED

 By Claims \ref{cl-wu} and \ref{cl-wvi}, we know that $N_{G^*}(w)\subseteq \{u_1,\ldots,u_a\}$ for any $w\in W$. Now, we can see that $G^*$ is also  $F_5$-free, and $G^*\neq S_{\frac{m+4}{2},2}^{-1}$ since $W\neq \varnothing$. By Theorem \ref{le-Chen}, $\lambda<\lambda(S_{\frac{m+4}{2},2}^{-1})$, a contradiction.  Hence, $W=\varnothing$.\QED
 
 \bigskip
 
 {\bf Concluding remarks.}   In this note, we characterize the extremal graph for $F_6$. For $k\geq 3$, the extremal graph for $F_{2K+1}$ or  $F_{2K+2}$ 
  only exists when $\frac{m}{k}-\frac{k-1}{2}$ is an integer. It is interesting to determine  extremal graphs for $F_{2K+1}$ or  $F_{2K+2}$  for the remaining case.

\end{document}